\input amstex
\documentstyle{amsppt}
\magnification=1200
\hyphenation{}
 \loadbold  
 \loadeufm  
\pagewidth{5.6in}
\pageheight{7.5in}
\catcode`\@=\active
\catcode`\@=11
\font\boldss=cmssbx10
\newfam\bssfam
\textfont\bssfam=\boldss
\def\bss{\relaxnext@\ifmmode\let\next\bss@\else
 \def\next{\Err@{Use \string\bss\space only in math mode}}\fi\next}
\def\bss@#1{{\bss@@{#1}}}
\def\bss@@#1{\fam\bssfam#1}
\catcode`\@=\active
\NoBlackBoxes
\topmatter
\title
A COMPACTNESS CRITERION \\
IN SCHAUDER SPACES
\endtitle
\leftheadtext{ M.  Lanza de Cristoforis}
\rightheadtext{{\rm A COMPACTNESS CRITERION IN SCHAUDER SPACES
}}
\author
MASSIMO LANZA DE CRISTOFORIS\footnote"\dag"{
Dipartimento di Matematica Pura ed Applicata, Universit\`a di Padova,
Via Belzoni 7, 35131 Padova, Italy.\newline}
\endauthor
\thanks
I am grateful to Prof\. Giuseppe De Marco who has
pointed me out the validity of Thorem 2.3 and has 
thus enabled me to write this note.
\endthanks 
\thanks
The present note had been pre-printed in 1991.
\endthanks
\abstract A necessary and
sufficient compactness criterion in Schauder Spaces is
proved.   
\endabstract
\endtopmatter
\document
\subhead{1. Introduction}\endsubhead
Although necessary and sufficient compactness criteria
in the space  $C^m(\text{cl}\,
\varOmega )$ of m-times continuously differentiable
functions on the closure $\text{cl}\,
\varOmega$ of the open subset $\varOmega$
of $\Bbb R^n$ are well-known,
(cf\. e\.g\. Kufner, John \& Fu\v cik (1977, p\. 37)),
it seems that only sufficient conditions for
compactness are known in the subspace $C^{m,\alpha}
(\text{\rm cl}\,\varOmega )$ of those $g\in C^m
(\text{\rm cl}\,\varOmega )$ whose derivatives of
order $m$ are H\" older continuous with exponent
$\alpha\in (0,1]$ 
(cf\. e\.g\. Kufner, John \& Fu\v cik (1977, p\.
38)).  By observing that a function $g\in C^{0,\alpha}
(\text{\rm cl}\,\varOmega )$ can be identified (up to
an additive constant) with the continuous and bounded
function $\left(g(\bold x)-g(\bold
y)\right)\cdot |\bold x-\bold y|^{-\alpha}$ on $(\bold
x,\bold y)\in \left(\varOmega\times\varOmega\right)
\setminus\left\{(\bold
x,\bold x):\bold x\in\varOmega\right\}$, we
characterize the relatively compact subsets of 
$C^{0,\alpha} (\text{\rm cl}\,\varOmega )$ by means of
the  Ascoli-Arzel\`a Theorem in the space of
continuous and bounded functions on
$\left(\varOmega\times\varOmega\right)\setminus\left\{(\bold
x,\bold x):\bold x\in\varOmega\right\}$. The criterion
of compactness is then easily extended to the space 
$C^{m,\alpha} (\text{\rm cl}\,\varOmega )$, $m\geq 1$.
\subsubhead{Notation} \endsubsubhead Let $\Cal
X$ be a Banach space. A subset $\Cal A\subseteq\Cal X$
is said to be relatively compact provided that its
closure $\text{\rm cl}\,\Cal A$ is compact. The
diameter of $\Cal A$ is denoted $\text{diam}\,\Cal A$.
 Let $\bold x_{0}\in\Cal X$, $r>0$. $B(\bold
x_{0},r)$ denotes $\left\{\bold x\in\Cal X:
\|\bold x-\bold x_{0}\|<r\right\}$. We denote by $A^E$
the set of  all functions from the set $E$ to the
set $A$. Let $C\subseteq A$ and  $f\in A^E$. We
denote $\left\{\bold x\in E:f(\bold x)\in C\right\}$
by $f^{\leftarrow }(C)$. We set $\Delta\equiv
\left\{(\bold x,\bold x):\,\bold x\in\Bbb R^n\right\}$. Let
$\varOmega$ be an open subset of $\Bbb R^n$. Let
 $C^0_b(\varOmega )$ denote the (Banach) space of all
real-valued continuous and bounded functions on
$\varOmega$, equipped with the supremum norm
$\|f\|_0\equiv\sup_{\bold x\in\varOmega }|f(\bold x
)|$. The  space of $m$-times continuously
differentiable functions on $\varOmega $, is denoted 
$C^m(\varOmega )$. Let 
$\eta\equiv (\eta_1,\dots ,\eta_n)\in\Bbb N^n, |\eta
|\equiv \eta_1+\dots +\eta_n  $. Then
$D^{\eta} $ denotes
$\frac{\partial^{|\eta|}}{\partial
x_1^{\eta_1}\dots\partial x_n^{\eta_n}}$.  The subspace
of $C^m(\varOmega )$ of those functions which have a continuous extension 
to $\text{\rm cl}\,\varOmega $ together with their
derivatives $D^{\eta }f$ of order $|\eta |\leq m$ is
denoted $C^m(\text{\rm cl}\,\varOmega )$. Let $f\in
C^m(\text{\rm cl}\,\varOmega )$. The unique continuous
extension of $D^\eta f$, $|\eta |\leq m$ to $
\text{\rm cl}\,\varOmega$ is still denoted by the same
symbol. Let $\varOmega $ be a bounded open subset of
$\Bbb R^n$.  $C^m(\text{\rm cl}\,\varOmega )$ equipped
with the norm $\|f\|_m\equiv \sum _{|\eta |\leq
m}\sup_{\text{\rm cl}\,\varOmega}|D^\eta f|$ is a
Banach space.   The subspace of $C^m(\text{\rm
cl}\,\varOmega ) $  whose functions have $m$-th order
derivatives that are H\"older continuous  with
exponent  $\alpha\in (0,1]$ is denoted $C^{m,\alpha}
(\text{\rm cl}\,\varOmega )$. If $f$ is a map of
$\text{cl}\, \varOmega $ to $\Bbb R^n$, then its 
H\"older quotient is $|f|_{\alpha }\equiv \sup\left\{
\frac{|f(\bold x )-f(\bold y)|}{|\bold x-\bold y|^\alpha
}:\bold x,\bold y\in \text{cl}\,\varOmega , \bold x\neq
\bold y\right\}$. The space $C^{m,\alpha }(\text{cl}\,
\varOmega )$ is equipped with its usual norm 
$\|f\|_{m,\alpha }=
\|f\|_m+\sum_{|\bold\eta |=m}
|D^{\bold\eta }f |_{\alpha}$.

\subhead{2. The Compactness Criterion}\endsubhead
Let $\varOmega$ be a bounded open subset of $\Bbb R^n$
and let $S$ be the map of $ C^{0,\alpha }(\text{cl}\,
\varOmega )$
into $C^0_b(\left(\varOmega\times\varOmega\right)\setminus\Delta )$ 
defined by 
$$
S[g](\bold x,\bold y)\equiv\frac{g(\bold x)-g(\bold
y)}{|\bold x-\bold y|^\alpha},\qquad g\in C^{0,\alpha
}(\text{cl}\, \varOmega ).
\tag 2.1
$$
The map $S$ is a linear isometry of the
space $ C^{0,\alpha }(\text{cl}\,
\varOmega )$ equipped with the seminorm
$|\cdot |_\alpha$to the space
$C^0_b\left(\left(\varOmega\times\varOmega\right)\setminus\Delta 
\right)$ equipped with the supremum norm
$\|\cdot\|_0$. Now, let $\overline{\bold
x}\in\text{cl}\,
\varOmega$. Since the norm $\|\cdot\|_{0,\alpha }$ is
clearly equivalent to the norm $\|g\|_{0,\alpha
}^{\overline{\bold x}}\equiv|g(\overline{\bold
x})|+|g|_{\alpha }$, the following holds.
\proclaim{Lemma 2.2} Let $\varOmega$ be a bounded 
open subset of $\Bbb R^n$
and let $\Cal H\subseteq  C^{0,\alpha
}(\text{cl}\, \varOmega )$. Then the following
conditions are equivalent
\roster
\item"{(i)}"There exists $\overline{\bold
x}\in\text{cl}\,\varOmega$ such that $\sup_{h\in\Cal
H}|h(\overline{\bold
x})|<+\infty$ and $S[\Cal H]$ is relatively
compact in $C^0_b\left(\left(\varOmega\times\varOmega\right)\setminus\Delta 
\right)$. 
\item"{(ii)}"The set $\Cal H$ is relatively compact
in  $C^{0,\alpha }(\text{cl}\,
\varOmega )$.
\endroster
\endproclaim
The relatively compact subsets of  
$C^0_b\left(\left(\varOmega\times\varOmega\right)\setminus\Delta 
\right)$ can be characterized by means of the
following variant of the Ascoli-Arzel\`a Theorem.
\proclaim{Theorem 2.3} Let $W$ be a subset of
$\Bbb R^n$, $\Cal K$ a subset of $C^0_b(W)$. The set
$\Cal K$ is relatively compact with respect to the
topology induced by the supremum norm $\|\cdot\|_0$
if and only if both the following conditions are
satisfied.
\roster
\item"{(i)}" For all $\bold x\in W$, the set
$\left\{f(\bold x): f\in\Cal K\right\}$
is relatively compact in $\Bbb R$.
\item"{(ii)}" For every $\epsilon >0$, there exists
a finite covering $\Cal F$ of subsets $F$ of $W$,
with $F$ open in the topology induced on $W$ by $\Bbb
R^n$, such that
$\text{\rm diam} f(F)<\epsilon$, $\forall F\in\Cal
F$, $\forall f\in\Cal K$.
\endroster
\endproclaim
\demo{Proof} This Theorem could be deduced by an
application of the Ascoli-Arzel\`a Theorem to the
space of continuous functions defined on the Stone-\v
Cech compactification of $W$ 
(cf\. Engelking (1977, 8.3.18), Isbell (1964, p\.
51).) However, to make the presentation elementary and
self-containing, we present a direct proof, which is
modelled on the classic proof of the  Ascoli-Arzel\`a
Theorem in the space $C^0(W)$, with 
$W$ compact, (cf\. e\.g\. 
Folland (1984, p\. 131).) Assume that $\Cal K$ is
relatively compact. Then $(i)$ clearly holds and $\Cal
K$ is totally bounded. Let $\epsilon >0$ and
$\{f_1,\dots ,f_r\}$ be an $\frac{\epsilon}{3}$-net
for $\Cal K$. Namely, let $\Cal
K\subseteq\cup_{i=1}^r\left\{f\in C^0_b(W):
\|f-f_i\|_0<\frac{\epsilon}{3}\right\}$. Now let $J$
be a compact subset of $\Bbb R$ such that
$J\supseteq\cup_{i=1}^rf_i(W)$ and let $\cup_{j=1}^s
B(y_j,\frac{\epsilon}{6})$ be a finite covering of open
balls in $\Bbb R$ of $J$. Then
$\cup_{j=1}^sf_i^{\leftarrow}\left(
B(y_j,\frac{\epsilon}{6})\right)$ is, for each $i$, a
finite open covering of $W$. Now let
$\left\{F_k\right\}_{k=1}^t$ be the finite covering
of $W$ obtained by taking the finite nonempty
intersections of $\left\{f_i^{\leftarrow}\left(
B(y_j,\frac{\epsilon}{6})\right) \right\}_{\Sb
i=1,\dots,r \\j=1,\dots,s\endSb}$. Clearly for each
$i$, each $F_k$ is contained in 
$f_i^{\leftarrow}\left(
B(y_j,\frac{\epsilon}{6})\right)$ for some $j$. Now 
let $\epsilon>0$. For all $f\in\Cal K$ there exists
$f_i$ such that $\|f-f_i\|_0<\frac{\epsilon}{3}$.
Then for each $f\in \Cal K$, $\bold x$, $\bold y\in
F_k$, we have 
$|f(\bold x)-f(\bold y)|\leq\allowmathbreak |f(\bold
x)-f_i(\bold x)|\allowmathbreak +|f_i(\bold
x)-f_i(\bold y)|+|f_i(\bold y)-f(\bold
y)|\allowmathbreak<\epsilon$.
\par Conversely, we now show that if conditions
$(i)$, $(ii)$ are satisfied, then $\Cal K$ is totally
bounded. Let $\epsilon>0$, and let $\Cal
F\equiv\left\{F_l\right\}_{l=1,\dots,\Lambda}$ be a
finite covering of $W$ as in $(ii)$ for
$\frac{\epsilon}{3}$. Let $\bold x_l\in F_l$,
$l=1,\dots ,\Lambda$, and 
$S\equiv\{\bold x_l : l=1,\dots ,\Lambda\}$. Since
for each $l$, the set $\left\{f(\bold x_l): f\in\Cal
K\right\}$ is relatively compact, then the set $\Cal
K_{|S}$ of the restrictions to $S$ of the functions
of $\Cal K$ is relatively compact in $\Bbb R^S$. Then
for all $\epsilon>0$ there exists a finite number of
elements $f_\theta\in\Cal K$, $\theta
=1,\dots,\Theta$ such that $\Cal
K\subseteq\cup_{\theta=1}^{\Theta}\left\{f\in\Cal
K:\sup_{l=1,\dots,\Lambda}|f(\bold
x_l)-f_{\theta}(\bold
x_l)|<\frac{\epsilon}{3}\right\}$. Now let $\bold
x\in W$ and $f\in \Cal K$. For some
$\overline{l}$ we have $\bold x\in F_{\overline{l}}$,
and for some $\theta$, $|f(\bold
x_{\overline{l}})-f_{\theta}(\bold
x_{\overline{l}})|<\frac{\epsilon}{3}$. Since by
assumption, $|f(\bold
x)-f(\bold
x_{\overline{l}})|<\frac{\epsilon}{3}$, and
$|f_{\theta}(\bold x)-f_{\theta}(\bold
x_{\overline{l}})|<\frac{\epsilon}{3}$, we conclude
that $|f_{\theta}(\bold
x)-f(\bold x)|<\epsilon$. By the 
arbitrariness of $\bold x$, we have $\Cal
K\subseteq\cup_{\theta=1}^{\Theta}\left\{f\in\Cal
K:\|f-f_{\theta}\|_0<\epsilon\right\}$.\qed
\enddemo 
Then, by combining Lemma 2.2  and Theorem 2.3, we have
the following. 
\proclaim{Theorem 2.4} Let $\varOmega$
be a bounded open subset of $\Bbb R^n$, $\alpha\in
(0,1]$. Let $\Cal K\subseteq C^{0,\alpha }(\text{cl}\,
\varOmega )$. Then the following conditions are
equivalent.
\roster
\item"{(i)}" The set $\Cal K$ is relatively compact
in $\left(C^{0,\alpha }(\text{cl}\,
\varOmega ),\|\cdot\|_{0,\alpha }\right)$.
\item"{(ii)}" For all $\bold
x\in\varOmega$, the set
$\left\{k(\bold x): k\in\Cal K\right\}$ is
relatively compact in $\Bbb R$, and for all $\epsilon
>0$, there exists a finite covering $\Cal F$ of open
subsets $F$ of
$\left(\varOmega\times\varOmega\right)\setminus\Delta$
such that
$\text{{\rm diam}}\left(S[k](F)\right)<\epsilon$,
$\forall F\in\Cal F$, $\forall k\in\Cal K$.  
\endroster
\endproclaim 
If $\varOmega $ is connected, then for
all $\bold x$, $\bold y\in\varOmega$ there exists an
arc $\gamma _{\bold x,\bold y} $ of class $C^1$ such
that   $\gamma _{\bold x,\bold y }:[0,1]\rightarrow
\varOmega$, $\gamma _{\bold x,\bold y }(0)=\bold x$,
$\gamma _{\bold x,\bold y }(1)=\bold y$. The geodesic
distance  $\lambda (\bold x,\bold y )$ is defined as 
$\lambda (\bold x,\bold y )\equiv  \inf \{\text{length
of }\gamma _{\bold x,\bold y }: \gamma _{\bold x,\bold
y } \text{is of class}\   C^1,\ \gamma
_{\bold x,\bold y }(0)=\bold x,\ \gamma_{\bold
x,\bold y }(1)=\bold y\}$. Let $c[\varOmega ]\equiv
\sup \left\{
\lambda (\bold x,\bold y )\cdot|\bold x
-\bold y |^{-1}:\bold x, \bold y\in\varOmega ,\ \bold
x\neq\bold y\right\}$. 
By using the inequality  $|f|_1\leq c[\varOmega ]
\sum_{i=1}^{n}\|\frac{\partial f}{\partial
x_i}\|_0$,  a simple inductive argument, and the
classic  Ascoli-Arzel\`a Theorem, 
we easily see that if
$\varOmega $ is a bounded open connected subset of
$\Bbb R^n$ such that $c[\varOmega ]<+\infty$, $m\geq
1$, then  $C^{m}(\text{cl}\,\varOmega )$ is
compactly imbedded in $C^{m-1}(\text{cl}\,\varOmega
)$. Then by the well-known
 compactness of the imbedding $C^{0,\alpha }
(\text{cl}\,\varOmega )\subseteq 
C^0(\text{cl}\,\varOmega )$, (which holds
also if $c[\varOmega ]=+\infty $) we deduce that  if 
$c[\varOmega ]<+\infty$ the space 
$C^{m,\alpha }(\text{cl}\,\varOmega )$ is compactly
imbedded in  $C^m(\text{cl}\,\varOmega )$, $m\geq 1$,
(which is also known.) Then we have the following. 
\proclaim{Theorem 2.5} Let $\varOmega$ be an
open, bounded connected subset of $\Bbb R^n$ such that
$c[\varOmega ]<+\infty $, $\alpha\in (0,1]$, $m\in\Bbb
N\setminus\{0\}$. Let $\Cal K\subseteq C^{m,\alpha
}(\text{cl}\, \varOmega )$. Then the following
conditions are equivalent. \roster
\item"{(i)}" The set $\Cal K$ is relatively compact
in $\left(C^{m,\alpha }(\text{cl}\,
\varOmega ),\|\cdot\|_{m,\alpha }\right)$.
\item"{(ii)}" The set $\Cal K$ is bounded in $\left(C^{m,\alpha }(\text{cl}\,
\varOmega ),\|\cdot\|_{m,\alpha }\right)$, and for
all $\epsilon >0$ there exists a finite covering
$\Cal F$ of open subsets $F$ of
$\left(\varOmega\times\varOmega\right)\setminus\Delta$
such that
$\text{{\rm diam}}\left(S[D^{\eta
}k](F)\right)\allowmathbreak<\epsilon$,  $\forall
F\in\Cal F$, $\forall k\in\Cal K$, $\forall\eta\in\Bbb
N^n$ with $|\eta|=m$.  \endroster 
\endproclaim
\demo{Proof} If $\Cal K$ is relatively compact in
$C^{m,\alpha }(\text{cl}\,\varOmega)$, then it is
obviously bounded and $\left\{D^{\eta}k:k\in\Cal K
\right.$, $\left.|\eta |=m\right\}$ is relatively
compact in $C^{0,\alpha }(\text{cl}\, \varOmega )$.
Then condition $(ii)$ holds by Theorem 2.4.
Conversely, let $\{k_n\}$ be a sequence in $\Cal K$.
By $(ii)$ and Theorem 2.4, for every $\eta\in\Bbb
N^n$, $|\eta |=m$ there exist a subsequence
$\left\{D^{\eta }k_{n_j}\right\}$ and an element
$k^{\eta }$ of $C^{0,\alpha }(\text{cl}\, \varOmega )$
such that $\lim_j D^{\eta }k_{n_j}=k^{\eta }$ in
$C^{0,\alpha }(\text{cl}\,\varOmega)$. Since
$c[\varOmega]<+\infty$, the space  $C^{m,\alpha
}(\text{cl}\,\varOmega )$ is compactly embedded in
$C^m(\text{cl}\, \varOmega )$. Then by possibly
considering a subsequence, we can assume that there
exists $k\in C^m(\text{cl}\,\varOmega )$ such that
$\lim_jk_{n_j}=k$ in  $C^m(\text{cl}\,\varOmega)$.
Then $D^{\eta }k=k^{\eta }$ and 
 $\lim_jk_{n_j}=k$ in 
$C^{m,\alpha }(\text{cl}\,\varOmega)$.\qed\enddemo

\enddocument